\documentclass[12pt]{article}
\usepackage{amssymb}
\usepackage{amsfonts}
\usepackage{amsmath}
\usepackage[usenames]{color}
\usepackage{mathrsfs}
\usepackage{amsfonts}
\usepackage{amssymb,amsmath}
\usepackage{CJK}
\usepackage{cite}
\usepackage{cases}
\usepackage{amsthm}

\def\cl{\centerline}

\def\vs{\vspace*}

\def\ni{\noindent}

\numberwithin{equation}{section}
\newtheorem{theo}{Theorem}[section]
\newtheorem{defi}[theo]{Definition}

\newtheorem{lemm}[theo]{Lemma}
\newtheorem{exam}[theo]{Example}

\newtheorem{remark}[theo]{Remark}

\begin{document}
\begin{center}
\cl{\large\bf \vs{6pt} Isoparametric hypersurfaces and hypersurfaces }
\cl{\large\bf \vs{6pt} with constant principal curvatures in Finsler spaces}
\footnote {$^*\,$ Project supported by AHNSF (No.2108085MA11).
\\\indent\ \ $^\dag\,$ chenylwuhu@qq.com
}
\cl{Peilong Dong$^1$, Yali Chen$^2$$^\dag\,$ }

\cl{\small 1. School of Mathematics and Statistics, Zhengzhou Normal University, } \cl{\small Zhengzhou Henan 450044, China.}
\cl{\small 2. School of Mathematics and Statistics, Anhui Normal University, }  \cl{\small Wuhu, Anhui, 241000, China. }
\end{center}

{\small
\parskip .005 truein
\baselineskip 3pt \lineskip 3pt

\noindent{{\bf Abstract:} In this paper, we study the relationship between isoparametric hypersurfaces and hypersurfaces
 with constant principal curvatures in Finsler spaces. We give some examples of isoparametric hypersurfaces with (non)constant principal curvatures on Randers manifolds with nonconstant flag curvatures. Furthermore, we construct an example of a conformally flat Randers manifold which admits a family of nonisoparametric hyperplanes with constant principal curvatures.

\vs{5pt}

\ni{\bf Key words}: isoparametric hypersurfaces, constant principal curvatures, mean curvatures, parallel hypersurfaces, nonconstant flag curvatures}

\ni{\it Mathematics Subject Classification (2010)}: 53C60, 53B40. }
\parskip .001 truein\baselineskip 6pt \lineskip 6pt
\section{Introduction}

In Riemannian geometry, a hypersurface is called isoparametric if itself and all
its sufficiently close parallel hypersurfaces have constant mean
curvature. Since 1938, E. Cartan began to study the isoparametric hypersurfaces of real space forms systematically and
proved that a hypersurface has constant principal curvatures if and only if it is isoparametric \cite{C38}.

In Finsler geometry, the conception of isoparametric hypersurfaces has been firstly introduced in \cite{HYS16}. Moreover, the analytic definitions of
locally isoparametric hypersurfaces
and $d\mu$-isoparametric hypersurfaces are given in~\cite{HDY, DH}.  Let~$(N,F,d\mu)$ be an
$n$-dimensional Finsler manifold with volume form~$d\mu=\sigma(x)dx$ and $f$ be a non-constant $C^1$ function and smooth on
$N_f=\{x\in N~|~df(x)\neq 0\}$.
Set $J=f(N_f)$. Then $f$ is
said to be \emph{isoparametric} (resp. \emph{$d\mu$-isoparametric}) on $(N,F,d\mu)$,
if there exist a smooth function ${a} (t)$ and a continuous function $b (t)$ on $J$ such that
\begin{equation}\label{1.3} \left\{\begin{aligned}
&F(\nabla f)=a(f),\\
&\Delta f=b(f),
\end{aligned}\right.
\end{equation}
where $\nabla f$ denotes the gradient of~$f$ and $\Delta {f}=\Delta^{\nabla f} f$ (resp. $\Delta {f}=\Delta_{\sigma} f$) is
defined by (\ref{l1})(resp. (\ref{l2})).
All the regular level surfaces ${M}_t = {f}^{-1}(t)$ are named an \emph{($d\mu$-)isoparametric family}, each of which is called an \emph{($d\mu$-)isoparametric hypersurface} on $(N,F,d\mu)$. A function $f$ only satisfying the first equation
of (1.1) is said to be \emph{transnormal}.

Similar to Riemannian geometry, \cite{HYS16} proved that
a transnormal function $f$ is isoparametric if and only if each regular level hypersurface of $f$ has constant ($d\mu$)-mean curvatures. Particularly, if $N$ has constant flag curvature (and constant $\mathbf{S}$-curvature), then a transnormal function $f$ is ($d\mu$)-isoparametric if and only if all the principal curvatures of $M_{t}$ are constant. So
the relationship between isoparametric hypersurfaces and principal curvatures is worth studying.

The previous work mainly focused on isoparametric hypersurfaces with constant principal curvatures in Finsler space forms (with constant flag curvature). For some very special Finsler space forms, such as Minkowski spaces (with zero flag curvature) and Funk spaces (with negative constant flag curvature), the isoparametric hypersurfaces have been completely classified in \cite{HYS16, GM, HYS17}. In~\cite{HDRY}, the authors studied isoparametric hypersurfaces of Finsler space forms, obtained the Cartan-type formula and some classifications on the number of distinct principal curvatures or their multiplicities.

In Riemannian geometry, the equivalence between isoparametric hypersurfaces and hypersurfaces with constant principal curvatures is no longer true for more general ambient spaces with nonconstant curvatures. There exist many examples of isoparametric hypersurfaces with (non)constant principal
curvatures with nonconstant curvatures \cite{RC13, GTY}. In \cite{R19}, a conformally flat metric in $\mathbb{R}^{n}$ is constructed, which admits a family of nonisoparametric hyperplanes. This shows that there is no necessary relation between the isoparametric hypersurface and the hypersurface with constant principal curvature.

A natural question is whether there exist similar results in Finsler spaces. In this paper, we give some examples of isoparametric hypersurfaces with (non)constant principal curvatures in Randers manifolds with nonconstant flag curvatures. Furthermore, we construct an example of a conformally flat Finsler manifold which admits a family of hypersurface with constant principal curvatures, but each of them is an nonisoparametric hypersurface, that is, its
sufficiently close parallel hypersurfaces don't have constant mean curvatures or $d\mu_{BH}$-mean curvatures. So we have the following theorem.
\begin{theo}
For Finsler spaces with nonconstant flag curvatures, the relation between isoparametric hypersurfaces and hypersurfaces with constant principal curvatures is no longer equivalent. In fact, there exist isoparametric hypersurfaces without constant principal curvatures and nonisoparametric hypersurfaces with constant principal curvatures.
\end{theo}
\section{Preliminaries}

\subsection{Finsler manifolds}
Let~$(N,F)$ be an $n$-dimensional smooth connected Finsler manifold and~$TN$ be the tangent bundle over~$N$ with local coordinates
$(x,y)$, where~$x=(x^i)$ and~$y=y^i\frac{\partial}{\partial x^{i}}$.
Here and from now on, we will use the following convention of index ranges unless otherwise stated:
$$1\leq i, j, \cdots \leq n ;~~~~~~~1\leq a, b, \cdots \leq n-1.$$
The fundamental form~$g$ of~$(N,F)$ is
\begin{equation*}
g=g_{ij}(x,y)dx^{i} \otimes dx^{j}, ~~~~~~~g_{ij}(x,y)=\frac{1}{2}[F^{2}] _{y^{i}y^{j}}.
\end{equation*}
The projection~$\pi : TN\rightarrow N$ gives rise to the pull-back bundle~$\pi^{\ast}TN$. On it there exists a unique
\emph{Chern connection}~$\nabla$ with~$\nabla
\frac{\partial}{\partial x^i}=\omega_{i}^{j}\frac{\partial}{\partial
x^{j}}=\Gamma^{i}_{jk}dx^k\otimes\frac{\partial}{\partial x^{j}}$ satisfying (\cite {BCS})
$$dg_{ij}-g_{ik}\omega^{k}_{j}-g_{kj}\omega^{k}_{i}=2C_{ijk}(dy^{k}+N^{k}_{l}dx^{l}),$$
$$~~~~N_j^i:=\frac{\partial G^i}{\partial y^j}=\Gamma^{i}_{jk}y^k,$$
where $C_{ijk}=\frac{1}{2}\frac{\partial g_{ij}}{\partial y^k}$ is
called the\emph{ Cartan tensor} and $G^{i}=\frac{1}{4}g^{il}\left\{[F^{2}]_{x^{k}y^{l}}y^{k}-[F^{2}]_{x^{l}}\right\}$ are the \emph{geodesic coefficients} of $(N,F)$. For $X=X^{i}\frac{\partial}{\partial x^{i}}\in\Gamma(TN)$, the covariant derivative of $X$ along $v=v^{i}\frac{\partial}{\partial x^{i}}\in T_{x}N$ with respect to a reference vector $w\in T_{x}N\setminus\{0\}$ is defined by
$$D^{w}_{v}X(x)=\{v^{j}\frac{\partial X^{i}}{\partial x^{j}}(x)+\Gamma^{i}_{jk}(w)v^{j}X^{k}(x)\}\frac{\partial}{\partial x^{i}}.$$
Another torsion-free Berwald connection $^{b}\nabla$ is defined by
$$^{b}\omega^{i}_{j}=^{b}\Gamma^{i}_{jk}dx^{k}=(\Gamma^{i}_{jk}+L^{i}_{jk})dx^{k}$$
where $L^{i}_{jk}$ is the Landsberg curvature of $(N, F)$. For~$X=X^{i}\frac{\partial}{\partial x^{i}}\in \pi^{\ast}TN$, the \emph{horizontal covariant derivative} of~$X$ with respect to $^{b}\nabla$ is defined by
\begin{align}
X^{i}_{~|j}=\frac{\delta X^{i}}{\delta x^{j}}+^b\Gamma^{i}_{js}X^s,\label{Z10}
\end{align}
where $\frac{\delta}{\delta x^{j}}=\frac{\partial}{\partial x^{j}}-N_j^i\frac{\partial}{\partial y^{i}}$.

Let~${\mathcal L}:TN\rightarrow T^{\ast}N$ denotes the \emph{Legendre transformation}, satisfying~${\mathcal L}(\lambda
y)=\lambda {\mathcal L}(y)$ for all~$\lambda>0,~y\in TN$.
For a smooth function~$f: N\rightarrow \mathbb{R}$, the \emph{gradient vector} of~$f$ at~$x$ is defined as~$\nabla f(x)={\mathcal
L}^{-1}(df(x))\in T_{x}N$.
Set~$N_{f}=\{x\in N|df(x)\neq 0\}$ and~$\nabla^{2}f(x)=\nabla^{\nabla f}(\nabla f)(x)$ for~$x\in N_{f}$. We define the \emph{Laplacian} of~$f$
by
\begin{equation}\label{l1}
\Delta^{\nabla f} f=\textmd{tr}_{g_{_{\nabla f}}}(\nabla^{2}f).
\end{equation}
And the \emph{Laplacian} of~$f$ with respect to the volume form~$d\mu=\sigma(x)dx=\sigma(x)dx^{1}\wedge dx^{2}\wedge\cdots\wedge dx^{n}$ can be represented as
\begin{align}\label{l2}
\Delta_{\sigma} f=\textmd{div}_{\sigma}(\nabla f)=\frac{1}{\sigma}\frac{\partial}{\partial x^{i}}(\sigma g^{ij}(\nabla f)f_{j})=\Delta^{\nabla f} f-S(\nabla f),
\end{align}
where
\begin{align}\label{1.4.19}
S(x,y)=\frac{\partial G^{i}}{\partial y^{i}}-y^{i}\frac{\partial}{\partial x^{i}}(\ln \sigma(x))
\end{align}
is the \emph{$\mathbf{S}$-curvature}.

Let $\phi:M\to N$ be an embedded hypersurface of~$(N,F)$. For any~$x\in M$,
there exist exactly two unit\emph{ normal vectors}~$\mathbf{n}_{\pm}$. Let~$\mathbf{n}$ be a given normal vector of $N$.
Set $\hat g=\phi^*g_{\mathbf{n}}$.
From \cite{HYS16}, we have the following \emph{Gauss-Weingarten formulas}
\begin{align}\label{1.1}
{ D}^{\textbf{n}}_{X}Y&={ {\hat{\nabla}}}_{X}Y+\hat{h}(X,Y)\textbf{n},\\
D^{\textbf{n}}_{X}\textbf{n}&=-{A}_{\textbf{n}}X,~~~~~\quad \forall X,~Y\in \Gamma(TM),\label{1.2}
\end{align}
where $\hat{\nabla}$ is the induced connection on $(M,\hat{g})$ and $\hat{h}(X,Y):=g_{\textbf{n}}(\textbf{n},D^{\textbf{n}}_{X}Y)$. The eigenvalues of the shape operator ${A}_{\mathbf{n}}$,
$\mu_1,\mu_2,\cdots,\mu_{n-1}$, and~$\hat{H}_{\mathbf{n}}=\sum\limits_{a=1}^{n-1}\mu_{a}$ are called the \emph{ principal
curvatures} and the \emph{anisotropic mean curvature} with respect to~$\mathbf{n}$, respectively.

\begin{lemm}\cite{HYS16} Let $M$ be an embedded hypersurface of $(N, F, d\mu)$, then
\begin{equation}\label{hn}
H_{\textbf{$\mathbf{n}$}}=\hat{H}_{\textbf{$\mathbf{n}$}}+S(\textmd{\textbf{$\mathbf{n}$}}),
\end{equation}
where $H_{\textbf{$\mathbf{n}$}}$ is the $d\mu$-mean curvature induced by variation of the volume in $(N, F, d\mu)$ (see (3.3) in \cite{HYS16} for detail).
\end{lemm}
Similar to Riemannian geometry, we can also give the geometric definitions of isoparametric hypersurface on Finsler manifolds as following.
\begin{defi}\label{d1}
A hypersurface $M$ on Finsler manifold $(N, F)$ is called isoparametric if it and its sufficiently close parallel hypersurfaces have constant anisotropic mean curvatures.
Moreover,
a hypersurface $M$ in Finsler manifold $(N, F, d\mu)$ is named $d\mu$-isoparametric if it and all its sufficiently close parallel hypersurfaces have constant $d\mu$-mean curvatures.
\end{defi}
\begin{remark}
Obviously, when $S(\textmd{\textbf{$\mathbf{n}$}})=\textmd{constant}$, a hypersurface $M$ is $d\mu$-isoparametric if and only if it is isoparametric.
\end{remark}

\section{The examples of isoparametric hypersurfaces on Randers manifolds with nonconstant flag curvatures}

In Riemannian geometry, there are many examples of isoparametric hypersurfaces with (non)constant principal
curvatures on some Riemannian manifolds with nonconstant curvature \cite{RC13, GTY, W82, DRDV12}.
The first examples were found by Wang \cite{W82}, who constructed some inhomogeneous isoparametric
hypersurfaces with nonconstant principal curvatures in the complex projective space, by projecting some of the
inhomogeneous hypersurfaces in spheres via the Hopf map. In \cite{DRDV12}, the examples of isoparametric real hypersurfaces with nonconstant principal curvatures in complex hyperbolic spaces $\mathbb{C}H^{n}$ are constructed by Lie groups and Lie algebras.

In this section, we plan to give some examples of isoparametric hypersurfaces with (non)constant principal curvatures in Randers manifolds with nonconstant flag curvatures.

In section 2, we study the principal curvatures of anisotropic submanifolds and hypersurfaces in a Randers
space with the navigation datum. Here we will use this relationship to construct examples.
Let~$(N, h)$ be an~$n$-dimensional Riemannian manifold and $v$ be a vector
field. By navigation problem, $(h, v)$ can define a Randers metric
\begin{equation}\label{c3}
F=\frac{\sqrt{\lambda h^2+v_{0}^2}-v_{0}}{\lambda}=\frac{\sqrt{\lambda h_{ij}y^iy^j+(v_{i}y^i)^2}-v_{i}y^i}{\lambda},
\end{equation}
where~$\lambda=1-b^2, ~b=\|v\|_h, ~ v=v^i\frac{\partial}{\partial x^{i}}, ~v_i=h_{ij}v^j$.
Let~$\phi: M\to(N, F)$ be an embedded hypersurface and the unit normal vector field of~$M$ with respect to~$F$ and $h$ be $\mathbf{n}$ and $\bar{\mathbf{n}}$, respectively. From \cite{HDY}, we know
\begin{equation}
\mathbf{n}=\bar{\mathbf{n}}+v\label{a3.21}
\end{equation}
and the relationship between the principal curvatures of~$M$ in~$(N, F)$ and $(N, h)$ is as follows.
\begin{lemm} \label{thm0}\cite{HDY}
Let~$M$ be a hypersurface in a Randers space~$(N,F,d\mu_{BH})$ with the navigation datum~$(h, v)$. If~$F$ has isotropic~$\mathbf{S}$-curvature~$S=(n+1)k(x)F$, then for any unit normal vector field $\mathbf{n}$, the shape operators of $M$ in Randers space~$(N, F)$ and Riemannian space~$(N, h)$, ${A}_{\mathbf{n}}$ and $\bar{A}_{\bar{\mathbf{n}}}$,  have the same principal vectors and their principal curvatures satisfy
\begin{align}\label{3.24}\mu=\bar{\mu}+k(x),\end{align}
where~$\mu$ and~$\bar{\mu}$ are the principal curvatures of~$M$ in Randers space~$(N, F)$ and Riemannian space~$(N, h)$, respectively.
\end{lemm}
Later, Xu and his coworkers got the local correspondence between isoparametric functions
or isoparametric hypersurfaces for homothetic navigation in a Finsler manifold.
\begin{lemm}\cite{XMYZ}\label{x1}
Let $v$ be a homothetic vector field on Finsler manifold $(N, \bar{F})$, and $F$ the metric
defined by navigation from the datum $(\bar{F}, v )$. Then
locally around any point $x_{0}$ with $\bar{F}(x_{0}, -v(x_{0}))< 1$, a hypersurface is isoparametric
for $(\bar{F}, d\mu^{\bar{F}}_{BH})$ if and only if it is isoparametric for $(F, d\mu^{F}_{BH})$.
\end{lemm}

According to \cite{H05}, Randers metric $F$  with respect to the Busemann-Hausdorff (BH) measure has constant $\mathbf{S}$-curvature if and only if $v$ is a homothetic vector field.
From Lemma \ref{thm0} and Lemma \ref{x1}, we know that if $v$ is a homothetic vector field, the ($d\mu^{F}_{BH}$-)isoparametric hypersurface of $(N, h)$ is also ($d\mu^{F}_{BH}$-)isoparametric
for the Randers manifold $(N, F)$, and the principal curvatures of $M$ in $(N, F)$ are all (non)constant
if and only if its principal curvatures in $(N, h)$ are all (non)constant.

\begin{remark}\label{k1}
By using homothetic navigation, we can construct ($d\mu^{F}_{BH}$-)isoparametric hypersurfaces with (non)constant principal
curvatures on some Randers manifolds with nonconstant flag curvatures.
\end{remark}
From Lemma \ref{thm0}, Lemma \ref{x1} and \cite{DRDV12}, we have the following example.
\begin{exam}
Let $\mathbb{C}H^{n}$ be a complex hyperbolic space, $AN$ be the connected and simply connected
subgroups of $G=SU(1, n)$, $S_{\varpi}$ be a
connected subgroup of $AN$, $M^{r}$ be
the tube of radius $r$ around the submanifold $W_{\varpi}$, where $W_{\varpi}$ is the orbit $S_{\varpi}\cdot o$ of the
group $S_{\varpi}$ through the point $o\in\mathbb{C}H^{n}$. Then for every $r > 0$, $M^{r}$
is an isoparametric real hypersurface which has, in general, nonconstant principal curvatures. If
given a (locally) homothetic vector field $v$ in $\mathbb{C}H^{n}$, then by (\ref{c3}), we can get a Randers metric $F$ with nonconstant flag curvatures.
So $M^{r}$ is also
the isoparametric hypersurface of $(\mathbb{C}H^{n}, F)$, and the principal curvatures of $M^{r}$ in $(\mathbb{C}H^{n}, F)$ are all nonconstant.
\end{exam}

In Riemannian spaces of constant curvatures, a hypersurface is isoparametric if and only if it has constant principal curvatures.
But in Finsler manifolds with constant flag curvatures,
for any given volume form $d\mu$ on $(N, F)$, if $F$ doesn't have constant $\mathbf{S}$-curvature, then
the existence of isoparametric hypersurfaces with nonconstant principal curvatures is still an open problem.
\section{The example of nonisoparametric hypersurfaces with constant principal curvatures}
Let~$N=\mathbb{R}^{n}$, $h=\sqrt{\delta_{ij}y^{i}y^{j}}$ and $v=b\frac{\partial}{\partial x^{n}}$
be a vector field, where $b$ is a constant and satisfies $b<1$.
Consider a Randers-Minkowski metric
\begin{equation}
F(y)=\alpha+\beta=\frac{\sqrt{\lambda h^2+v^2_{0}}}{\lambda}-\frac{v_{0}}{\lambda},
\end{equation}
where
$ v_{0}=by^n, \lambda=1-\delta_{ij}v^{i}v^{j}=1-b^2$.
Set $\alpha=\sqrt{\alpha_{ij}y^{i}y^{j}}, ~\beta=b_{i}y^{i}$, then
\begin{equation*}
\alpha_{ab}=\frac{\delta_{ab}}{\lambda}, ~\alpha_{an}=0, ~\alpha_{nn}=\frac{\lambda+b^{2}}{\lambda^2}=\frac{1}{\lambda^2},
\end{equation*}
\begin{equation*}
\alpha^{ab}=\lambda\delta_{ab}, ~\alpha^{an}=0, ~\alpha^{nn}=\lambda^2,
\end{equation*}
\begin{equation*}
b^{a}=\alpha^{aj}b_{j}=0, ~b^{n}=-\lambda b, ~\parallel\beta\parallel^{2}_\alpha=b^{2}.
\end{equation*}
Hence
$$\alpha_{y^{a}}=\frac{y^{a}}{\lambda\alpha},~\alpha_{y^{n}}=\frac{y^{n}}{\lambda^{2}\alpha},$$
\begin{equation}\label{gab}
g_{ab}=\frac{F}{\lambda\alpha}\delta_{ab}-\frac{\beta y^{a}y^{b}}{\alpha^{3}\lambda^{2}}, ~g_{an}=-\frac{by^{a}}{\lambda^{2}\alpha}-\frac{\beta y^{a}y^{n}}{\lambda^{3}\alpha^{3}},
\end{equation}
\begin{equation*}
g_{nn}=\frac{F}{\lambda^{2}\alpha}\left(1-\frac{(y^{n})^{2}}{\lambda^{2}\alpha^{2}}\right)+(\frac{y^{n}}{\lambda^{2}\alpha}-\frac{b}{\lambda})^{2},
\end{equation*}
and
\begin{equation*}
g^{ab}=\frac{\alpha}{F}\lambda\delta^{ab}+\frac{b^{2}\alpha+\beta}{F^{3}}y^{a}y^{b},
~g^{an}=\frac{b\lambda\alpha}{F^{2}}y^{a}+\frac{b^{2}\alpha+\beta}{F^{3}}y^{a}y^{n},
\end{equation*}
\begin{equation*}
g^{nn}=\frac{\lambda^{2}\alpha}{F}+2\frac{b\lambda\alpha}{F^{2}}y^{n}+\frac{b^{2}\alpha+\beta}{F^{3}}(y^{n})^{2}.
\end{equation*}
\subsection{A conformally flat Randers manifold}
Define a new metric
\begin{equation*}\label{dl1}
\widetilde{F}(x, y)=e^{\rho(x)}F(y)
\end{equation*}
and the corresponding fundamental tensor is
\begin{equation}
\tilde{g}_{ij}(x, y)=e^{2\rho(x)}g_{ij}(y),
\end{equation}
where $\rho(x)=\sum\limits_{a=1}^{n-1} \ln(2+cos(\pi x^{a}))$, $x=(x^{i})\in \mathbb{R}^{n}$.
It is easy to know that $\tilde{g}$ is (locally) conformally flat.
From \cite{SS16}, we get the relation of spray coefficient $\tilde{G}^{i}$ and $G^{i}$ with respect to $\widetilde{F}$ and $F$, respectively,
\begin{equation}\label{gg}
\tilde{G}^{i}=G^{i}+\rho_{k}y^{k}y^{i}-\frac{1}{2}F^{2}\rho^{i},
\end{equation}
where $\rho_{k}=\frac{\partial \rho}{\partial x^{k}}, \rho^{i}=g^{ij}\rho_{j}$.
By $\rho_{n}=0$, $G^{i}=0$,
\begin{equation*}
\rho^{a}=g^{ai}\rho_{i}=g^{ab}\rho_{b}=\frac{\alpha}{F}\lambda\rho_{a}+\frac{ b^{2}\alpha+\beta}{F^{3}}y^{a}y^{b}\rho_{b}
\end{equation*}
and
\begin{equation*}
\rho^{n}=g^{nb}\rho_{b}=(\frac{b\lambda\alpha}{F^{2}}y^{b}+\frac{b^{2}\alpha+\beta}{F^{3}}y^{b}y^{n})\rho_{b}.
\end{equation*}
We have
\begin{equation}\label{g1}
\tilde{G}^{a}=\rho_{b}y^{b}y^{a}-\frac{1}{2}F^{2}\rho^{a}
=(1-\frac{b^{2}\alpha+\beta}{2F})\rho_{b}y^{b}y^{a}-\frac{1}{2}\lambda\alpha F \rho_{a},
\end{equation}
\begin{equation}\label{g2}
\tilde{G}^{n}=\rho_{b}y^{b}y^{n}-\frac{1}{2}F^{2}\rho^{n}=(1-\frac{b^{2}\alpha+\beta}{2F})\rho_{b}y^{b}y^{n}
-\frac{b\lambda\alpha}{2}\rho_{b}y^{b}.
\end{equation}
Then
by (\ref{gg}), it followed that
\begin{equation}\label{n1}
\tilde{N}^{i}_{j}=\frac{\partial \tilde{G}^i}{\partial y^j}
=\rho_{k}y^{k}\delta^i_j+\rho_{j}y^i
-\frac{1}{2}(F^{2})_{y^{j}}\rho^i+F^2 C^{ik}_{j}\rho_{k},
\end{equation}
where
\begin{align}\label{n2}
&-2C^{ik}_{j}\nonumber=\frac{\partial g^{ik}}{\partial y^{j}}
=(\frac{\alpha}{F})_{y^{j}}\alpha^{ik}-(\frac{\alpha}{F^2})_{y^{j}}(b^{i}y^{k}+b^{k}y^{i}),\\
&-\frac{\alpha}{F^2}(b^{i}\delta^{k}_{j}+b^{k}\delta^{i}_{j})
+(\frac{b^{2}\alpha+\beta}{F^3})_{y^{j}}y^{i}y^{k}+\frac{b^{2}\alpha+\beta}{F^3}(y^{i}\delta^{k}_{j}+y^{k}\delta^{i}_{j}).
\end{align}
From \cite{SS16}, we have
\begin{equation*}
\tilde{S}(y)=S(y)+F^{2}\rho^{k}I_{k},
\end{equation*}
where $I_{k}=g^{ij}C_{ijk}$ is the Cartan form, $\tilde{S}$ and $S$ respectively are the $\mathbf{S}$-curvature of $\tilde{F}$ and $F$.
It is well known that $S=0$ with respect to the Busemann-Hausdorff volume form.
Then by (\ref{n2}), $\rho_{n}=0$ and $\rho_{k}b^{k}=0$, we have
\begin{align}\label{s1}
\tilde{S}\nonumber=&F^{2}\rho_{k}C^{ik}_{i}\\ \nonumber
=&-\frac{1}{2}F^{2}\rho_{k}\left[(\frac{\alpha}{F})_{y^{i}}\alpha^{ik}
-(\frac{\alpha}{F^2})_{y^{i}}b^{i}y^{k}
+(\frac{b^{2}\alpha+\beta}{F^3})_{y^{i}}y^{i}y^{k}+(n+1)y^{k}\frac{b^{2}\alpha+\beta}{F^3}\right]\\
=&-\frac{1}{2}\rho_{k}\left[(\alpha_{y^{i}}\beta-\alpha\beta_{y^{i}})\alpha^{ik}-b^{i}y^{k}\frac{\alpha_{y^{i}}(\beta-\alpha)-2\alpha\beta_{y^{i}}}{F}
+(n-1)y^{k}\frac{b^{2}\alpha+\beta}{F}\right].
\end{align}

\subsection{The geodesics with respect to $\widetilde{F}$}
By (\ref{dl1}), $\tilde{F}$ is produced by navigation datum $(\tilde{h}, \tilde{v})$,
where $\tilde{h}=e^{\rho}h, \tilde{v}=e^{\rho}v$. Then
$\tilde{v}_{a}=0, \tilde{v}_{n}=b e^{\rho}$, $\tilde{v}^{a}=0, \tilde{v}^{n}=b e^{-\rho}$.
Let $M=\{x\in\mathbb{R}^{n}|x^{n}=x_{0}^{n}\in\mathbb{R}\}$, $\bar{\mathbf{n}}$ and $\tilde{\mathbf{n}}$ be the unit normal vector field of~$M$ with respect to~$\tilde{h}$ and $\tilde{F}$, respectively.
According to \cite{HDY}, we have
\begin{equation*}\label{3.21}
\tilde{\mathbf{n}}=\bar{\mathbf{n}}+\tilde{v}.
\end{equation*}
From \cite{R19}, we know that $\bar{\mathbf{n}}=(0, 0, \cdots, e^{-\rho})$, then $\tilde{\mathbf{n}}=(1+b)\bar{\mathbf{n}}$.

Define
\begin{equation*}
\Omega:=\{(m_{1}, m_{2}, \cdots, m_{n-1}, x^{n})\in N|~m_{a}\in\mathbb{Z}\}.
\end{equation*}
Let $p=(m_{1}, m_{2}, \cdots, m_{n-1}, x_{0}^{n})\in\Omega$, $\{e_{1}, e_{2}, \cdots, e_{n-1}, \tilde{\textbf{n}}\}$ be a local orthogonal frame in the neighbourhood of $p$ with respect to $\tilde{F}$. Let $\gamma$ be an unit-speed geodesic with respect to $\tilde{F}$ and start at $p$ with initial direction $\tilde{\textbf{n}}$. By (\ref{g1}) and (\ref{g2}), we know that
$\tilde{F}$ is not projectively flat. So the geodesics with respect to $\tilde{F}$ are not necessarily straight lines.
Set

$(1)~\Lambda_{a}: x=(x^{1}, \cdots, x^{a}, \cdots, x^{n})\rightarrow \hat{x}=(x^{1}, \cdots, -x^{a}, \cdots, x^{n})$,

$(2)~\Psi_{m_{a}}: x=(x^{1}, \cdots, x^{a}, \cdots, x^{n})\rightarrow \check{x}=(x^{1}, \cdots, x^{a}+2m_{a}, \cdots, x^{n})$.

Next we will prove $\hat{\gamma}(t)=\Lambda_{a}(\gamma(t))$
and $\check{\gamma}(t)=\Psi_{m_{a}}(\gamma(t))$ are all geodesics.
So $\hat{\gamma}(t)$ and $\check{\gamma}(t)$ should satisfy
\begin{equation}\label{c1}
\ddot{\gamma}^{i}+\widetilde{\Gamma}^{i}_{jk}\dot{\gamma}^{k}\dot{\gamma}^{j}=0.
\end{equation}
For $\Lambda_{a}$, we get
$$\ddot{\hat{\gamma}}^{a}=-\ddot{\gamma}^{a}, ~\alpha(\dot{\hat{\gamma}})=\alpha(\dot{\gamma}),~F(\dot{\hat{\gamma}})=F(\dot{\gamma}),$$
$$\rho_{a}(\hat{\gamma})=-\rho_{a}(\gamma),~\rho_{b}(\hat{\gamma})\dot{\hat{\gamma}}^{b}=\rho_{b}(\gamma)\dot{\gamma}^{b}.$$
So by (\ref{g1}), we have
\begin{align*}
&\ddot{\hat{\gamma}}^{b}+2\widetilde{G}^{b}(\hat{\gamma}, \dot{\hat{\gamma}})\\ \nonumber
&=\ddot{\hat{\gamma}}^{b}+2\left[\left(1-\frac{b^{2}\alpha(\dot{\hat{\gamma}})+\beta(\dot{\hat{\gamma}})}{2F( \dot{\hat{\gamma}})}\right)\rho_{c}(\hat{\gamma})\dot{\hat{\gamma}}^{c}\dot{\hat{\gamma}}^{b}-\frac{1}{2}\lambda F(\dot{\hat{\gamma}})\alpha(\dot{\hat{\gamma}})\rho_{b}(\hat{\gamma})\right]\\ \nonumber
&=\left\{\begin{aligned}
      &-(\ddot{\gamma}^{b}+2\widetilde{G}^{b}(\gamma, \dot{\gamma})), ~a=b,\\
      &\ddot{\gamma}^{b}+2\widetilde{G}^{b}(\gamma, \dot{\gamma}), ~~~~~~~a\neq b.
\end{aligned}\right.
\end{align*}
Then $\hat{\gamma}(t)=\Lambda_{a}(\gamma(t))$ is a geodesic.

For $\Psi_{m_{a}}$, we also get
$$\ddot{\check{\gamma}}^{a}=\ddot{\gamma}^{a},
~F(\dot{\check{\gamma}})=F(\dot{\gamma}), ~\rho_{a}(\check{\gamma})=\rho_{a}(\gamma).$$
Then
\begin{align*}
&\ddot{\check{\gamma}}^{a}+2\widetilde{G}^{a}(\check{\gamma}, \dot{\check{\gamma}})\\ \nonumber
&=\ddot{\check{\gamma}}^{a}+2\left[\left(1-\frac{b^{2}\alpha(\dot{\check{\gamma}})+\beta(\dot{\check{\gamma}})}{2F( \dot{\check{\gamma}})}\right)\rho_{c}(\check{\gamma})\dot{\check{\gamma}}^{c}\dot{\check{\gamma}}^{a}-\frac{1}{2}\lambda\alpha(\dot{\check{\gamma}})F( \dot{\check{\gamma}})\rho_{a}(\check{\gamma})\right]\\ \nonumber
&=\ddot{\gamma}^{a}+2\widetilde{G}^{a}(\gamma, \dot{\gamma}).
\end{align*}
So $\hat{\gamma}(t)=\Psi_{a}(\gamma(t))$ is also a geodesic.

From the above, $\tilde{\gamma}(t)=\Psi_{m_{a}} \Lambda_{a}(\gamma(t))
=(\gamma^{1}(t), \cdots, -\gamma^{a}(t)+2m_{a}, \cdots, \gamma^{n}(t))$
is also a geodesic with respect to $\tilde{F}$. But $\tilde{\gamma}(t)$ and $\gamma(t)$ have the same initial condition.
Hence, by the uniqueness of geodesic, we get $\gamma^{a}(t)=m_{a}\in\mathbb{Z}$.
Note that $\rho(\gamma(t))=c\ln 3$,
where $c$ is the number of even entries of $(m_{1}, m_{2}, \cdots, m_{n-1})$.
Since
\begin{equation}\label{rou1}
\rho_{a}(\gamma)=\frac{-\pi sin(\pi x^{a})}{2+cos(\pi x^{a})}(\gamma)=0.
\end{equation}
Then by (\ref{g1}) and (\ref{g2}), we get
$$\widetilde{G}^{a}({\gamma}, \dot{\gamma})=\widetilde{G}^{n}({\gamma}, \dot{\gamma})=0. $$
From (\ref{c1}), we have
\begin{equation}
\ddot{\gamma}^{n}=0.
\end{equation}
According to initial condition, then
\begin{equation}\label{j1}
\gamma(t)=(m_{1}, m_{2}, \cdots, m_{n-1}, x_{0}^{n}+(1+b)e^{-\rho}t),
\end{equation}
where $t$ is arc length.

\subsection{Parallel hypersurfaces}

For $x_{0}\in M$, there exist a neighbourhood $U(x_{0})\subset M$ and a distance function $r(x)$ to $M$, such that $r(x)=t$.
Then $\tilde{\mathbf{n}}=\nabla r$ and all the integral curves of $\nabla r$ are all normal geodesics (see \cite{HYS16} for detail).
So we can define the parallel hypersurfaces of $M$ by
$$M_{t}=\{x\in M~|~r(x)=t, t>0\}.$$

For $X=X^{i}\frac{\partial}{\partial x^{i}}\in TM$, note that
$X^{n}=0$ and $\tilde{F}_{y^{j}}X^{j}=F_{y^{j}}X^{j}=\alpha_{y^{j}}X^{j}=0$.
Then by $\alpha(\tilde{\mathbf{n}})=\frac{1}{1-b}e^{-\rho}, \beta(\tilde{\mathbf{n}})=\frac{-b}{1-b}e^{-\rho}$, (\ref{n1}) and (\ref{n2}), we get
\begin{align*}
D^{\tilde{\mathbf{n}}}_{X}\tilde{\mathbf{n}}
=&(\tilde{\mathbf{n}}^{i}_{x^{j}}+\tilde{N}^{i}_{j}(\tilde{\mathbf{n}}))X^{j}\frac{\partial}{\partial x^{i}}\nonumber\\
=&\left[\tilde{\mathbf{n}}^{i}_{x^{a}}+\rho_{a}\tilde{\mathbf{n}}^i
+\frac{1}{2}F^{2}(\tilde{\mathbf{n}})\rho_{a}\left(\frac{\alpha(\tilde{\mathbf{n}})}{F(\tilde{\mathbf{n}})^2}b^{i}
-\frac{b^{2}\alpha+\beta}{F(\tilde{\mathbf{n}})^3}\tilde{\mathbf{n}}^{i}\right)\right]X^{a}\frac{\partial}{\partial x^{i}} \nonumber\\
=&(\tilde{\mathbf{n}}^{i}_{x^{a}}+\rho_{a}\tilde{\mathbf{n}}^i)X^{a}\frac{\partial}{\partial x^{i}} \nonumber\\
=&0.
\end{align*}
By Weingarten formula (\ref{1.2}), $M$ has constant principal curvatures $\tilde{\theta}=0$ with respect to $\tilde{F}$.
From \cite{HYS16} (Lemma 4.4), we have
\begin{equation*}
\frac{d k_{a}(t)}{dt}=\tilde{K}(\tilde{\mathbf{n}}, e_{a})+k^{2}_{a}(t),
\end{equation*}
where $k_{a}(t)$ are the principal curvature of $M_{t}$, $k_{a}(0)=\mu_{a}$ and $\tilde{K}$ is the flag curvature of $(N, \tilde{F})$. Take the trace
\begin{equation}\label{h1}
\frac{d}{dt}\sum\limits_{a=1}^{n-1}k_{a}(t)=\widetilde{Ric}(\tilde{\mathbf{n}})+\sum\limits_{a=1}^{n-1}k_{a}^{2}(t),
\end{equation}
where $\widetilde{Ric}$ is the Ricci curvature of $(N, \tilde{F})$.

Now we prove that $M_{t}$ is neither isoparametric nor $d\mu$-isoparametric.
By \cite{SC07}, we have
\begin{align}
\widetilde{Ric}(y)
=&Ric(y)-F^{2}\{2\rho_{k}J^{k}+g^{ij}\rho_{i|j}+[\rho^{j}I_{j}]_{|k}y^k+2\theta\rho^{j}I_{j}\}\\ \nonumber
&+(n-2)(\Xi-F^{2}\parallel d\rho\parallel^{2}_{F})+F^{4}\rho^{j}\rho^{k}\{2I_{s}C^{s}_{ij}-[I_{j}]_{y^{k}}-C^{s}_{ij}C^{i}_{ks}\}, \nonumber
\end{align}
where $Ric(y)$ is the Ricci curvature of $(N, F)$, $J^{k}=g^{ij}L^{k}_{ij}, ~\theta(y)=\rho_{k}y^{k},$
$\Xi=\theta^{2}-\theta_{x^k}y^{k}$, $\parallel d\rho\parallel^{2}_{F}=g^{ij}\rho_{i}\rho_{j}.$
Note that $G^{i}=0$,
$L_{ijk}=0$, $Ric(y)=0$, $C_{ijn}(\tilde{\mathbf{n}})=0$. By (\ref{gab}), we get
\begin{equation*}
2C_{abc}=-\frac{(\delta_{ab}y^{c}+\delta_{ac}y^{b}+\delta_{bc}y^{a})\beta}{\lambda^{2}\alpha^3}
+3\frac{\beta y^{a}y^{b}y^c}{\lambda^3 \alpha^5},
\end{equation*}
\begin{equation*}
2C_{abn}=\frac{\delta_{ab}(\alpha\beta_{y^{n}}-\beta\alpha_{y^{n}})}{\lambda\alpha^2}
-\frac{y^{a}y^{b}(\alpha\beta_{y^{n}}-3\beta\alpha_{y^{n}})}{\lambda^2 \alpha^4},
\end{equation*}
then $C_{abc}(\tilde{\mathbf{n}})=0$. In conclusion, $C_{ijk}(\tilde{\mathbf{n}})=0$.
\begin{align*}
2\frac{\partial C_{abc}}{\partial y^{d}}
&=-\frac{(\delta_{ab}\delta_{cd}+\delta_{ac}\delta_{bd}+\delta_{bc}\delta_{ad})\beta}{\lambda^{2}\alpha^3}\\
&-3\frac{(\delta_{ab}y^{c}+\delta_{ac}y^{b}+\delta_{bc}y^{a})\beta\alpha_{y^{d}}}{\lambda^{2}\alpha^4}
+3(\frac{\beta y^{a}y^{b}y^c}{\lambda^3 \alpha^5})_{y^{d}},
\end{align*}
\begin{align*}
2\frac{\partial C_{abn}}{\partial y^{n}}
&=\frac{\delta_{ab}[-\beta\alpha\alpha_{y^{n} y^{n}}-2(\alpha\beta_{y^{n}}-\beta\alpha_{y^{n}})\alpha_{y^{n}}]}{\lambda\alpha^3}
-y^{a}y^{b}[\frac{(\alpha\beta_y^{n}-3\beta\alpha_y^{n})}{\lambda^2 \alpha^4}]_{y^{n}}.
\end{align*}
Then $C_{abc}(\tilde{\mathbf{n}})=0, ~2\frac{\partial C_{abc}}{\partial y^{d}}(\tilde{\mathbf{n}})=\frac{(\delta_{ab}\delta_{cd}+\delta_{ac}\delta_{bd}+\delta_{bc}\delta_{ad})be^{2\rho}}{(1+b)^{2}}$,
$\frac{\partial C_{abn}}{\partial y^{n}}(\tilde{\mathbf{n}})=0$.
$$\theta(\tilde{\mathbf{n}})=0,
~F(\tilde{\mathbf{n}})=e^{-\rho}, ~\Xi(\tilde{\mathbf{n}})=0,
~\frac{\alpha(\tilde{\mathbf{n}})}{F(\tilde{\mathbf{n}})}=\frac{1}{1-b},~I_{i}(\tilde{\mathbf{n}})=0,$$
$$g^{ab}(\tilde{\mathbf{n}})=(1+b)\delta^{ab},~g^{an}(\tilde{\mathbf{n}})=0,~g^{nn}(\tilde{\mathbf{n}})=(1+b)^{2},$$
$$\parallel d\rho\parallel^{2}_{F}(\tilde{\mathbf{n}})=g^{ij}(\tilde{\mathbf{n}})\rho_{i}\rho_{j}
=(1+b)\delta^{ab}\rho_{a}\rho_{b}=(1+b)\sum\limits_{a=1}^{n-1}(\rho_{a})^2,$$
$$g^{ij}(\tilde{\mathbf{n}})\rho_{i|j}
=-(1+b)\delta^{ab}\frac{\pi^{2}+2\pi^{2}cos(\pi x^{b})}{[2+cos(\pi x^{a})]^{2}}
=-(1+b)\sum\limits_{a=1}^{n-1}\frac{\pi^{2}+2\pi^{2}cos(\pi x^{a})}{[2+cos(\pi x^{a})]^{2}},$$
$$[\rho^{j}I_{j}]_{|k}y^k(\tilde{\mathbf{n}})=(\rho^{j}_{|k}I_{j}(\tilde{\mathbf{n}})+\rho^{j}I_{j|k}(\tilde{\mathbf{n}}))\tilde{\mathbf{n}}^k
=\rho^{j}g^{il}\frac{\partial C_{ilj}}{\partial x^k }\tilde{\mathbf{n}}^k=0,$$
\begin{align*}
&\rho^{j}\rho^{k}[I_{j}]_{y^{k}}(\tilde{\mathbf{n}})
=\rho^{j}\rho^{k}[-2C^{il}_{k}C_{ilj}+g^{il}\frac{\partial C_{ijl}}{\partial y^{k}})](\tilde{\mathbf{n}})\\
&=\rho^{a}\rho^{b}[g^{cd}\frac{\partial C_{bcd}}{\partial y^{a}}+g^{cn}\frac{\partial C_{abc}}{\partial y^{n}}+g^{nn}\frac{\partial C_{abn}}{\partial y^{n}}](\tilde{\mathbf{n}})\\ \nonumber
&=\rho^{a}\rho^{b}\frac{(\delta_{ab}\delta_{cd}+\delta_{ac}\delta_{bd}+\delta_{bc}\delta_{ad})be^{2\rho}\delta^{cd}}{(1+b)}
=\sum\limits_{a=1}^{n-1}(\rho^{a})^2\frac{(n+2)be^{2\rho}}{(1+b)}.
\end{align*}
So
\begin{align}
&\widetilde{Ric}(\tilde{\mathbf{n}})=e^{-2\rho}(1+b)\sum\limits_{a=1}^{n-1}\left[\frac{\pi^{2}+2\pi^{2}cos(\pi x^{a})}{[2+cos(\pi x^{a})]^{2}}
-(n-2)(\rho_{a})^2-(\rho^{a})^2\frac{(n+2)b}{(1+b)^2}\right].
\end{align}
If $u_{1}=(0, 0, \cdots, x_{0}^{n})\in M\cap\Omega$ and $u_{2}=(1, 1, \cdots, x_{0}^{n})\in M\cap\Omega$, we know
$$\widetilde{Ric}(\tilde{\mathbf{n}})|_{\gamma_{u_{1}}}=3^{(-2n+1)}(1+b)\pi^{2}(n-1)>0,$$
and
$$\widetilde{Ric}(\tilde{\mathbf{n}})|_{\gamma_{u_{2}}}=-(1+b)\pi^{2}(n-1)<0.$$
Meanwhile, when $t=0$, $k_{a}^{2}|_{\gamma_{u_{1}}(0)}=k_{a}^{2}|_{\gamma_{u_{2}}(0)}=0$. Then by (\ref{h1}),
$\frac{d}{dt}\sum\limits_{i=1}^{n}k_{a}|_{\gamma_{u_{1}}(0)}>0$ and $\frac{d}{dt}\sum\limits_{i=1}^{n}k_{a}|_{\gamma_{u_{2}}(0)}<0$.
So the mean curvatures of the parallel hypersurfaces $M_{t}$ is not constant.
Note that $\rho_{k}\tilde{\mathbf{n}}^{k}=0$ and by (\ref{s1}), we have
$$\tilde{S}(\tilde{\mathbf{n}})=0.$$
\begin{remark}
According to Definition \ref{d1}, $M$ with constant principal curvatures is neither isoparametric hypersurface nor $d\mu$-isoparametric hypersurface.
\end{remark}

Peilong Dong\\
School of Mathematics and Statistics, Zhengzhou Normal University,\\
Zhengzhou Henan 450044, China.\\
E-mail: dpl2021@163.com\\

Yali Chen\\
School of Mathematics and Statistics, Anhui Normal University, \\
Wuhu, Anhui, 241000, China.\\
E-mail: chenylwuhu@qq.com\\

\end{document}